\documentclass[twoside,12pt]{article}
\usepackage{amssymb,amsmath,bm,euscript}
\usepackage{graphicx,color,caption}

\usepackage[colorlinks=true]{hyperref}

\usepackage{algorithm,algpseudocode}
\usepackage{braket}

\newtheorem{proposition}{Proposition}[section]
\newtheorem{theorem}[proposition]{Theorem}
\newtheorem{lemma}[proposition]{Lemma}
\newtheorem{corollary}[proposition]{Corollary}
\newtheorem{definition}[proposition]{Definition}

\newcommand{\qed}{\hphantom{.}\hfill $\Box$\medbreak}

\def\II{\mathcal{I}_{m, n}}

\def\A{{\mathcal{A}}}

\def\B{{\mathcal{B}}}
\def\C{\mathcal{C}}

\def\aa{{\bf a}}
\def\ee{{\bf e}}
\def\x{{\bf x}}
\def\y{{\bf y}}

\def\z{{\bf z}}
\def\w{{\bf w}}
\def\s{{\bf s}}
\def\uu{{\bf u}}
\def\vv{{\bf v}}
\def\0{{\bf 0}}

\title{\bf{Biquadratic Tensors, Biquadratic Decomposition and Norms of Biquadratic Tensors}}
\author{ \hspace{1mm} Liqun Qi\thanks{Department of Applied
    Mathematics, The Hong Kong Polytechnic University, Hung Hom,
    Kowloon, Hong Kong; ({\tt liqun.qi@polyu.edu.hk}). This author's work was supported by the Hong Kong
    Research Grant Council (Grant No.  PolyU 15300715, 15301716 and 15300717). },
    \ \
    Shenglong Hu\thanks{Department of Mathematics, School of Science, Hangzhou Dianzi University, Hangzhou 310018 China; ({\tt shenglonghu@hdu.edu.cn}). This author's work was supported by NSFC (Grant No.  11771328).}
 \ and\
    Xinzhen Zhang\thanks{School of Mathematics, Tianjin University, Tianjin 300354 China; ({\tt xzzhang@tju.edu.cn}). This author's work was supported by NSFC (Grant No.  11871369). }}

\begin{document}
\date{\today}
\maketitle

\begin{abstract}
Biquadratic tensors play a central role in many areas of science.    Examples include elasticity tensor and Eshelby tensor in solid mechanics, and Riemann curvature tensor in relativity theory.  The singular values and spectral norm of a general third order tensor are the square roots of the M-eigenvalues and spectral norm of a biquadratic tensor. The tensor product operation is closed for biquadratic tensors.   All of these motivate us to study biquadratic tensors, biquadratic decomposition and norms of biquadratic tensors.   We show that the spectral norm and nuclear norm for a biquadratic tensor may be computed by using its biquadratic structure.   Then, either the number of variables is reduced, or the feasible region can be reduced.   We show constructively that for a biquadratic tensor, a biquadratic rank-one decomposition always exists, and show that the biquadratic rank of a biquadratic tensor is preserved under an independent biquadratic Tucker decomposition.   We present a lower bound and an upper bound of the nuclear norm of a biquadratic tensor.   Finally, we define invertible biquadratic tensors, and present a lower bound for the product of the nuclear norms of an invertible biquadratic tensor and its inverse, and a lower bound for the product of the nuclear norm of an invertible biquadratic tensor, and the spectral norm of its inverse.


\end{abstract}


\section{Introduction}

In this paper, unless otherwise stated, all the discussions will be carried out in the field of real numbers.

Suppose that $m$ and $n$ are positive integers.  Without loss of generality, we may assume that $m \le n$.

As in \cite{KB09}, we use $\circ$ to denote the operation of tensor outer product.   Then for $\x \in \Re^m$ and $\y \in \Re^n$, $\x \circ \y \circ \x \circ \y$ is a fourth order rank-one tensor in $\Re^{m \times n \times m \times n}$.  By the following definition, it is actually a biquadratic rank-one tensor.

\begin{definition}
Let $\Re^{m \times n \times m \times n}$ be the space of fourth order tensors of dimension $m \times n \times m \times n$.   Let $\A = (a_{i_1j_1i_2j_2}) \in \Re^{m \times n \times m \times n}$.  The tensor $\A$ is called biquadratic if for all $i_1, i_2 = 1, \cdots, m$ and $j_1, j_2 = 1, \cdots, n$, we have
$$a_{i_1j_1i_2j_2} = a_{i_2j_1i_1j_2} = a_{i_1j_2i_2j_1}.$$
The tensor $\A$ is called positive semi-definite if for any $\x \in \Re^m$ and $\y \in \Re^n$,
$$\langle \A, \x \circ \y \circ \x \circ \y \rangle \equiv \sum_{i_1, i_2 =1}^m \sum_{j_1, j_2 = 1}^n a_{i_1j_1i_2j_2}x_{i_1}y_{j_1}x_{i_2}y_{j_2} \ge 0.$$
The tensor $\A$ is called positive definite if for any $\x \in \Re^m, \x^\top \x = 1$ and $\y \in \Re^n, \y^\top \y = 1$,
$$\langle \A, \x \circ \y \circ \x \circ \y \rangle \equiv \sum_{i_1, i_2 =1}^m \sum_{j_1, j_2 = 1}^n a_{i_1j_1i_2j_2}x_{i_1}y_{j_1}x_{i_2}y_{j_2} > 0.$$
\end{definition}

Denote the set of all biquadratic tensors in $\Re^{m \times n \times m \times n}$ by $BQ(m, n)$.   Then $BQ(m, n)$ is a linear space.

Biquadratic tensors play a central role in many areas of science.    Examples include the elasticity tensor and the Eshelby tensor in solid mechanics, and the Riemann curvature tensor in relativity theory.    The elasticity tensor may be the most well-known tensor in solid mechanics and engineering \cite{Ny85}.    The Eshelby inclusion problem is one of the hottest topics in modern solid mechanics \cite{ZHHZ10}.  Furthermore, the Riemann curvature tensor is the backbone of Einstein's general relativity theory \cite{Ri06}.

Biquadratic tensors have very special structures.  The tensor product of two biquadratic tensors are still a biquadratic tensor.  This makes them very special.   Biquadratic tensors also have an M-eigenvalue structure.   An important problem in solid mechanics is if strong ellipticity condition holds or not \cite{KS75, SS83}.  In 2009, M-eigenvalues were introduced for the elastic tensor to characterize the strong ellipticity condition in \cite{QDH09}.   An algorithm for computing the largest M-eigenvalue was presented in \cite{WQZ09}.   The biquadratic optimization problem was studied in \cite{LNQY09}.  The M-eigenvalue structure was further extended to the Riemann curvature tensor \cite{XQW18}.    As the big data era arrived, the tensor completion problem came to the stage.  It was shown that the nuclear norm of tensors plays an important role in the tensor completion problem \cite{YZ16}.   A typical model in the tensor completion problem for higher order models is a general third order tensor \cite{Hu15, YZ16}.   The nuclear norm is the dual norm of the spectral norm \cite{FL17, Hu15, YZ16}.  The spectral norm of a tensor is its largest singular value.   In \cite{QH19}, it was shown that if we make contraction of a third order tensor with itself on one index, then we get a positive semi-definite biquadratic tensor.   A real number is a singular value of that third order tensor if and only if it is the square root of an M-eigenvalue of that positive semi-definite biquadratic tensor.   Thus, the spectral norm of that third order tensor is the square root of the spectral norm of that positive semi-definite biquadratic tensor.

All of these make biquadratic tensors a research interest.   In this paper, we study biquadratic tensors, biquadratic decomposition and norms of biquadratic tensors.   In the next section, we show that the spectral norm and nuclear norm for a biquadratic tensor may be computed by using its biquadratic structure.   Then, either the number of variables is reduced, or the feasible region can be reduced.   In Section 3, we show constructively that for a biquadratic tensor, a biquadratic rank-one decomposition always exists.   This gives an upper bound for the biquadratic rank of a biquadratic tensor.   In Section 4, we show that  the biquadratic rank of a biquadratic tensor is preserved under an independent biquadratic Tucker decomposition.   In Section 5, we present a lower bound and an upper bound of the nuclear norm of a biquadratic tensor.   In Section 6, we define invertible biquadratic tensors, and present a lower bound for the product of the nuclear norms of an invertible biquadratic tensor and its inverse, and a lower bound for the product of the nuclear norm of an invertible biquadratic tensor, and the spectral norm of its inverse.   Some final remarks are made in Section 7.

We use small letters $\lambda, x_i, u_i$, etc., to denote scalars, small bold letters $\x, \uu, \vv$, etc.,  to denote vectors, capital letters $A, B, C$, etc., to denote matrices, and calligraphic letters $\A, \B, \C$, etc., to denote tensors.

\section{Norms and M-Eigenvalues of Biquadratic Tensors}

For a vector $\uu = (u_1, \cdots, u_m)^\top$, we use $\|\uu\|_2$ to denote its 2-norm.   Thus,
$$\| \uu\|_2 := \sqrt{u_1^2 + \cdots + u_m^2}.$$

For a tensor $\A \in \Re^{m \times n \times m \times n}$, its spectral norm is defined as \cite{FL17, Hu15, JYZ17, YZ16}
\begin{equation} \label{n1}
\| \A \|_S := \max \left\{ \left| \langle \A, \x \circ \y \circ \uu \circ \vv \rangle \right| :  \|\x\|_2 = \|\y\|_2 = \|\uu\|_2 = \|\vv\|_2 = 1, \x, \uu \in \Re^m, \y, \vv  \in \Re^n \right\}.
\end{equation}

We have the following theorem.

\begin{theorem}  \label{t01}
Suppose that $\A \in BQ(m, n)$. Then
\begin{equation} \label{n3}
\| \A \|_S = \max \left\{ \left| \langle \A, \x \circ \y \circ \x \circ \y \rangle \right| :  \|\x\|_2 = \|\y\|_2 = 1, \x \in \Re^m, \y \in \Re^n \right\}.
\end{equation}
\end{theorem}
{\bf Proof} Suppose that the maximum of (\ref{n1}) is attained at $\bar \x, \bar \y, \bar \uu$ and $\bar \vv$.   Then
$$\| \A \|_S = \max \left\{ \left| \langle \A, \x \circ \bar \y \circ \uu \circ \bar \vv \rangle \right| :  \|\x\|_2 = \|\uu\|_2  = 1, \x, \uu \in \Re^m  \right\}.$$
Note that this is a homogeneous quadratic optimization.  Then there is a $\hat \x \in \Re^m$ such that $\| \hat \x\|_2 = 1$ and
$$\| \A \|_S = \left| \langle \A, \hat \x \circ \bar \y \circ \hat \x \circ \bar \vv \rangle \right|.$$
Then
$$\| \A \|_S = \max \left\{ \left| \langle \A, \hat \x \circ \y \circ \hat \x \circ \vv \rangle \right| :  \|\y\|_2 = \|\vv\|_2  = 1, \y, \vv \in \Re^n  \right\}.$$
Again, this is a homogeneous quadratic optimization.  Then there is a $\hat \y \in \Re^n$ such that $\|\hat \y\|_2 = 1$ and
$$\| \A \|_S = \left| \langle \A, \hat \x \circ \hat \y \circ \hat \x \circ \hat \y \rangle \right|.$$
This proves (\ref{n3}).
\qed

 In this way,  $\| \cdot \|_S$ also defines a norm in $BQ(m, n)$.

 Recall that the nuclear norm of $\A \in \Re^{m \times n \times m \times n}$ is defined as
\begin{equation} \label{n0}
\| \A \|_* = \inf \left\{ \sum_{j=1}^r |\lambda_j| : \A = \sum_{j=1}^r \lambda_j \x^{(j)} \circ \y^{(j)} \circ \uu^{(j)} \circ \vv^{(j)}, {\|\x^{(j)}\|_2 = \|\y^{(j)}\|_2 = \|\uu^{(j)}\|_2 = \|\vv^{(j)}\|_2 = 1, \atop \x^{(j)}, \uu^{(j)} \in \Re^m, \y^{(j)}, \vv^{(j)} \in \Re^n, r \in \mathbb{N}}
\right\}.
\end{equation}
By Corollary 5.4 of \cite{FL17}, we have
\begin{equation} \label{n01}
\| \A \|_* = \min \left\{ \sum_{j=1}^r |\lambda_j| : \A = \sum_{j=1}^r \lambda_j \x^{(j)} \circ \y^{(j)} \circ \x^{(j)} \circ \y^{(j)}, {\|\x^{(j)}\|_2 = \|\y^{(j)}\|_2 = 1, \atop \x^{(j)} \in \Re^m, \y^{(j)} \in \Re^n, r \in \mathbb{N}}
\right\}.
\end{equation}
 It can be calculated as \cite{FL17, Hu15, JYZ17, YZ16}
 \begin{equation} \label{n2}
\| \A \|_* := \max \left\{ \left| \langle \A, \B \rangle \right| :  \|\B \|_S = 1, \B \in \Re^{m \times n \times m \times n} \right\}.
\end{equation}
 For a biquadratic tensor, we have the following theorem.
 \begin{theorem}  \label{t02}
Suppose that $\A \in BQ(m, n)$. Then
\begin{equation} \label{n4}
\| \A \|_* = \max \left\{ \left| \langle \A, \B \rangle \right| :  \|\B \|_S = 1, \B \in BQ(m, n) \right\}.
\end{equation}
\end{theorem}
{\bf Proof} Without loss of generality, assume that $\A$ is nonzero.  Suppose that the maximum of (\ref{n2}) is attained at $\bar \B = (\bar b_{i_1j_1i_2j_2}) \in \Re^{m \times n \times m \times n}$ with $\| \bar \B \| = 1$.
Let $\hat \B = (\hat b_{i_1j_1i_2j_2})$ with
$$\hat b_{i_1j_1i_2j_2} = {1 \over 4} \left(\bar b_{i_1j_1i_2j_2} + \bar b_{i_2j_1i_1j_2} + \bar b_{i_1j_2i_2j_1} + \bar b_{i_2j_2i_1j_1}\right).$$
Then $\hat \B \in BQ(m, n)$, $\| \hat \B \|_S \le 1$, and
$$\| \A \|_* =  \left| \langle \A, \bar \B \rangle \right| = \left| \langle \A, \hat \B \rangle \right|.$$
Since $\A$ is not a zero tensor.  This implies that $\hat \B$ is also not a zero tensor.   Then $\| \hat \B \|_S \not = 0$.  Let
$$\tilde \B = {\hat \B \over \| \hat \B \|_S}.$$
We have $\tilde \B \in BQ(m, n)$, $\| \tilde \B \|_S = 1$, and
$$\left| \langle \A, \tilde \B \rangle \right| \ge  \left| \langle \A, \bar \B \rangle \right|.$$
Since $\bar \B$ is a maximizer of (\ref{n2}), we have
$$\left| \langle \A, \tilde \B \rangle \right| =  \left| \langle \A, \bar \B \rangle \right|.$$
This proves (\ref{n4}).
\qed

These two theorems show that we may compute the spectral norm and nuclear norm for a biquadratic tensor by using its biquadratic structure.   Then, either the number of variables is reduced, or the feasible region of the maximization problem can be reduced.   Furthermore, a biquadratic tensor has its own M-eigenvalue structure which is closely related to its spectral norm.

\begin{definition}
Suppose that $\A = (a_{i_1j_1i_2j_2}) \in BQ(m, n)$.   A real number $\lambda$ is called an M-eigenvalue of $\A$ if there are vectors  $\x = (x_1, \cdots, x_m)^\top \in \Re^m, \y = (y_1, \cdots, y_n)^\top \in \Re^n$ such that the following equations are satisfied:
For $i_1 = 1, \cdots, m$,
\begin{equation} \label{e5}
\sum_{i_2=1}^m \sum_{j_1, j_2=1}^n a_{i_1j_1i_2j_2}y_{j_1}x_{i_2}y_{j_2} = \lambda x_{i_1};
\end{equation}
For $j_1 = 1, \cdots, n$,
\begin{equation} \label{e6}
\sum_{i_1,i_2=1}^m\sum_{j_2=1}^n a_{i_1j_1i_2j_2}x_{i_1}x_{i_2}y_{j_2} = \lambda y_{j_1};
\end{equation}
and
\begin{equation} \label{e7}
\x^\top \x = \y^\top \y = 1.
\end{equation}
Then $\x$ and $\y$ are called the corresponding M-eigenvectors.
\end{definition}

\begin{theorem} \label{t1}
Suppose that $\A = (a_{i_1j_1i_2j_2}) \in BQ(m, n)$.  Then its M-eigenvalues always exist.  The spectral norm of $\A$ is equal to the largest absolute value of its M-eigenvalues.   Furthermore, $\A$ is positive semi-definite if and only if all of its M-eigenvalues are nonnegative; $\A$ is positive definite if and only if all of its M-eigenvalues are positive.  If $\A$ is positive semi-definite, then its spectral norm is equal to its largest M-eigenvalue.
\end{theorem}

This theorem was proved in \cite{QH19}.

For $m = n =3$, the elastic tensor in solid mechanics falls in the category of biquadratic tensors, with one additional symmetric properties between indices $i_1$ and $j_1$.   Then, the positive definiteness condition of $\A$ corresponds the strong ellipticity condition in solid mechanics.

\section{Biquadratic Rank-One Decomposition}

Let $\A \in BQ(m, n)$. Then $\A$ has a rank-one decomposition in the form
$$\A = \sum_{k=1}^r \x^{(k)} \circ \y^{(k)} \circ \s^{(k)} \circ \w^{(k)},$$
where $\x^{(k)}, \s^{(k)} \in \Re^m, \y^{(k)}, \w^{(k)} \in \Re^n$ for $k = 1, \cdots, r$.   The smallest $r$ for such a rank-one decomposition is called the rank of $\A$.

On the other hand, if we have
\begin{equation} \label{e4.5}
\A = \sum_{k=1}^r \x^{(k)} \circ \y^{(k)} \circ \x^{(k)} \circ \y^{(k)},
\end{equation}
where $\x^{(k)} \in \Re^m, \y^{(k)} \in \Re^n$ for $k = 1, \cdots, r$ for some positive integer $r$, then we say that $\A$ has a biquadratic rank-one decomposition.   The smallest for such a biquadratic rank-one decomposition is called the biquadratic rank of $\A$.   Denote it by $BR(\A)$.   The question is if such a biquadratic rank-one decomposition always exists.   We may following the approach in \cite{CGLM08} to show this by introducing biquadratic polynomials and using them as a tool for the proof.  Corollary 5.4 of \cite{FL17} also implies this.   Here, we give a constructive proof.   This also gives an upper bound of the biquadratic rank of a biquadratic tensor.

\begin{theorem} \label{t2}
For $\A = (a_{i_1j_1i_2j_2}) \in BQ(m, n)$, such a biquadratic rank-one decomposition always exists.  We also have
$$BR(\A) \le mn \min \left\{ {m(m+1) \over 2}, {n(n+1) \over 2}\right\}.$$
\end{theorem}
{\bf Proof}  For $\A = (a_{i_1j_1i_2j_2}) \in BQ(m, n)$, define a matrix
$$P = (p_{st}) \in \Re^{{m(m+1) \over 2} \times {n(n+1) \over 2}}$$
by
$$p_{st} = a_{i_1j_1i_2j_2}$$
for
$$s = {i_1(i_1-1) \over 2} + i_2$$
and
$$t = {j_1(j_1-1) \over 2} + j_2,$$
with $i_1 \ge i_2 \ge 1, j_1 \ge j_2 \ge 1$, $i_1 = 1, \cdots, m$ and $j_1 = 1, \cdots, n$.
Then $P$ has a singular value decomposition
$$M = \sum_{k=1}^q \sigma_k \uu^{(k)}\left(\vv^{(k)}\right)^\top,$$
where
$\uu^{(k)} \in \Re^{m(m+1) \over 2}$, $\|\uu^{(k)}\|_2=1$, $\vv^{(k)} \in \Re^{n(n+1) \over 2}$, $\|\vv^{(k)}\|_2=1$, for $k = 1, \cdots, q$, and
$$q \le \min \left\{ {m(m+1) \over 2}, {n(n+1) \over 2} \right\}.$$
For $\uu^{(k)} \in \Re^{m(m+1) \over 2}$, we may fold it to a symmetric matrix $U^{(k)} \in \Re^{m \times m}$.  Similarly, for $\vv^{(k)} \in \Re^{n(n+1) \over 2}$, we may fold it to a symmetric matrix $V^{(k)} \in \Re^{n \times n}$.   Suppose that $U^{(k)}$ has an eigenvalue decomposition
$$U^{(k)} = \sum_{l_u=1}^m \lambda_{l_u} \x^{(k, l_u)}\left(\x^{(k, l_u)}\right)^\top,$$
where $\x^{(k, l_u)} \in \Re^m$, $\|\x^{(k, l_u)}\|_2 = 1$, for $l_u=1, \cdots, m$, $k = 1, \cdots. q$.
Similarly, suppose that $V^{(k)}$ has an eigenvalue decomposition
$$V^{(k)} = \sum_{l_v=1}^m \mu_{l_v} \y^{(k, l_v)}\left(\y^{(k, l_v)}\right)^\top,$$
where $\y^{(k, l_v)} \in \Re^n$, $\|\y^{(k, l_v)}\|_2 = 1$, for $l_v=1, \cdots, n$, $k = 1, \cdots. q$.
Then we have
$$\A = \sum_{k=1}^q \sum_{l_u=1}^m \sum_{l_v=1}^n \sigma_k \lambda_{l_u}\mu_{l_v}\x^{(k, l_u)} \circ \y^{(k, l_v)} \circ \x^{(k, l_u)} \circ \y^{(k, l_v)}.$$
We have the conclusions.
\qed

Clearly, the biquadratic rank of a biquadratic tensor is always not less than its rank.   In which cases are these two ranks equal?   We do not go to further discussion on this in this paper.

Let $\A = (a_{i_1j_1i_2j_2}) \in BQ(m, n)$. Fix $j_1, i_2$ and $j_2$, then we have an $m$-vector $\aa_{\cdot j_1i_2j_2}$.  Denote by $A^{(1)}$ the $m \times mn^2$ matrix whose column vectors are such $m$-vectors for $i_2 = 1, \cdots, m$ and $j_1, j_2 = 1, \cdots, n$.  Then $A^{(1)}$ is the matrix flattening of $\A$ by the first index.   Here we do not specify the order of such column vectors in $A^{(1)}$ as this is not related.   Denote  the rank of $A^{(1)}$ by $R_1(\A)$.  We may define $R_2(\A)$, $R_3(\A)$ and $R_4(\A)$, respectively.  They are the Tucker ranks of $\A$ \cite{JYZ17, KB09}.   Then we have $R_1(\A) = R_3(\A)$ and $R_2(\A) = R_4(\A)$.  Hence, only $R_1(\A)$ and $R_2(\A)$ are independent.    We also have $R_1(\A) \le m$ and $R_2(\A) \le n$.

Suppose that $\A$ has a biquadratic rank-one decomposition as (\ref{e4.5}).
Denote $X$ as an $m\times r$ matrix, whose column vectors are $\x^{(1)}, \cdots, \x^{(r)}$, and $Y$ as an $n \times r$ matrix, whose column vectors are $\y^{(1)}, \cdots, \y^{(r)}$.  Then, as in \cite{KB09}, we may denote (\ref{e4.5}) as
\begin{equation} \label{e4.6}
\A = [[X, Y]]_{BQ}.
\end{equation}

Let $\A = (a_{i_1j_1i_2j_2}) \in BQ(m, n)$.   Suppose that $\A$ has a biquadratic rank-one decomposition (\ref{e4.6}).  Denote the ranks of $X$ and $Y$ by $R(X)$ and $R(Y)$ respectively.  Then we have
\begin{equation} \label{e4.7}
R(X) = R_1(\A), \ R(Y) = R_2(\A).
\end{equation}

\section{Biquadratic Tucker Decomposition}

We may also extend Tucker decomposition \cite{DDV00, JYZ17, KB09} to biquadratic Tucker decomposition.   Denote $\times_k$ as the mode-$k$ (matrix) product \cite{DDV00, JYZ17, KB09}.

\begin{definition}
Let $\A \in BQ(m, n)$. Suppose that there are $\B \in BQ(d_1, d_2)$, and $P \in \Re^{m \times d_1}$ and $Q \in \Re^{n \times d_2}$  such that
\begin{equation} \label{e5.8}
\A = \B \times_1 P \times_2 Q \times_3 P \times_4 Q := [[\B; P,Q]]_{BQ}.
\end{equation}
Then (\ref{e5.8}) is called a biquadratic Tucker decomposition of $\A$.   The tensor $\B$ is called a biquadratic Tucker core of $\A$.  The matrices $P$ and $Q$ are called the factor matrices of this decomposition.   A biquadratic Tucker decomposition is said to be independent if $P$ and $Q$ have full column rank.   A biquadratic Tucker decomposition is said to be orthonormal if $P$ and $Q$ have orthonormal columns.
\end{definition}

Note that if the biquadratic Tucker decomposition (\ref{e5.8}) is independent, then $d_1 \le m$ and $d_2 \le n$.

De Lathauwer, De Moor and Vandewalle \cite{DDV00} proposed an algorithm to compute Tucker decomposition (HOSVD) for a given tensor.   If we apply their algorithm to a biquadratic tensor,  since the first and the third matrix flattenings are the same, the second and the fourth flattenings are the same, we obtain an orthonormal biquadratic Tucker decomposition.

A biquadratic Tucker decomposition is a Tucker decomposition \cite{JYZ17, KB09}.  Thus, a biquadratic Tucker core has the properties of a Tucker core.  For example, the rank of a biquadratic Tucker core $\B$ is the same as the rank of $\A$, if the biquadratic Tucker decomposition is independent \cite{JYZ17}.   Similarly, the Tucker ranks will also be preserved by an independent biquadratic Tucker decomposition.  The problem is if some biquadratic properties, such as the biquadratic rank, and M-eigenvalues will be preserved or not.

We now prove the following theorem.

\begin{theorem} \label{t5.2}
Suppose that $\A \in BQ(m, n)$ has a biquadratic Tucker decomposition (\ref{e5.8}) and it is independent.  Then
$$BR(\A) = BR(\B).$$
\end{theorem}

We first prove a lemma.

\begin{lemma} \label{l1}
If the biquadratic Tucker decomposition (\ref{e5.8}) is independent, then there are $\hat P \in \Re^{d_1 \times m}$ and $\hat Q \in \Re^{d_2 \times n}$  such that
\begin{equation} \label{e5.9}
\B = \A \times_1 \hat P \times_2 \hat Q \times_3 \hat P \times_4 \hat Q := [[\A; \hat P,\hat Q]]_{BQ}.
\end{equation}
\end{lemma}
{\bf Proof} Let $\hat P = (P^\top P)^{-1}P^\top$ and $\hat Q = (Q^\top Q)^{-1}Q^\top$.  The conclusion follows.
\qed

{\bf Proof of Theorem \ref{t5.2}.} Suppose that $\B$ has a biquadratic rank-one decomposition
\begin{equation} \label{e5.10}
\B = \sum_{k=1}^r \hat \x^{(k)} \circ \hat \y^{(k)} \circ \hat \x^{(k)} \circ \hat \y^{(k)},
\end{equation}
where $\hat \x^{(k)} \in \Re^{d_1}$, $\hat \y^{(k)} \in \Re^{d_2}$ for $k = 1, \cdots, r$.
Then $\A$ has a biquadratic rank-one decomposition
(\ref{e4.5}) with
$$\x^{(k)} = P\hat \x^{(k)}, \ \y^{(k)} = Q\hat \y^{(k)},$$
for $k = 1, \cdots, r$.
This shows that
$$BR(\A) \le BR(\B).$$
Since the biquadratic Tucker decomposition (\ref{e5.8}) is independent, by Lemma \ref{l1}, we have (\ref{e5.9}).   Thus, if $\A$ has a biquadratic rank-one decomposition
(\ref{e4.5}), then $\B$ has a biquadratic rank-one decomposition
(\ref{e5.10}), with
$$\hat \x^{(k)} = \hat P \x^{(k)}, \ \hat \y^{(k)} = \hat Q \y^{(k)},$$
for $k = 1, \cdots, r$.
This shows that
$$BR(\A) \ge BR(\B).$$
Hence, we have
$$BR(\A) = BR(\B).$$
\qed

Jiang, Yang and Zhang \cite{JYZ17} proved the following theorem (Theorem 7 of \cite{JYZ17}).

\begin{theorem} \label{t5.4}
Suppose that $\A \in BQ(m, n)$ has a biquadratic Tucker decomposition (\ref{e5.8}) and it is orthonormal.  Then an M-eigenvalue of $\B$ is an M-eigenvalue of $\A$, and a
nonzero M-eigenvalue of $\A$ is also an M-eigenvalue of $\A$.
\end{theorem}

By Theorem \ref{t1}, this shows that the spectral norm is preserved under an orthonormal biquadratic Tucker decomposition.

Thus, biquadratic Tucker decomposition has better properties.   It only involves two factor matrices $P$ and $Q$.  This makes it much simple.

\section{Lower and Upper Bounds of the Nuclear Norm of a Biquadratic Tensor}

Let $\A = (a_{i_1j_1i_2j_2}) \in \Re^{m \times n \times m \times n}$.   If we regard $i_1j_1$ as an index from $1$ to $mn$, and regard $i_2j_2$ as another index from $1$ to $mn$, then we have a matrix flattening $M = M(\A) \in \Re^{mn \times mn}$.
 Then there is a one to one relation between $\A \in \Re^{m \times n \times m \times n}$ and $M \in \Re^{mn \times mn}$.   Hence, we may also write $\A = \A(M)$ for $M \in \Re^{mn \times mn}$.     If $M \in \Re^{mn \times mn}$ is diagonal, then we also say that $\A = \A(M)$ is diagonal.    In particular, if $M$ is the identity matrix $I_{mn} \in \Re^{mn \times mn}$, then we denote $\A(I_{mn})$ as $\II$ and call it the identity tensor in $\Re^{m \times n \times m \times n}$.

In the other words, for a fourth order tensor $\A = (a_{i_1j_1i_2j_2}) \in \Re^{m \times n \times m \times n}$, an entry $a_{i_1j_1i_2j_2}$ is called a diagonal entry if $i_1=i_2$ and $j_1=j_2$.   Otherwise, it is called an off-diagonal entry.  Then a diagonal tensor in $\Re^{m \times n \times m \times n}$ is a biquadratic tensor in $BQ(m, n)$ such that all of its off-diagonal entries are $0$, while the identity tensor $\II \in \Re^{m \times n \times m \times n}$ is the diagonal biquadratic tensor in $BQ(m, n)$ such that all of its diagonal entries are $1$.

Denote the Frobenius norm of a fourth order tensor $\A \in \Re^{m \times n \times m \times n}$ by $\|\A \|_2$, and the Frobenius norm  of $M \in \Re^{mn \times mn}$ by $\|M\|_2$.  For $\A, \B \in \Re^{m \times n \times m \times n}$, we use
  $\langle M(\A), M(\B) \rangle$ to denote the inner product of matrices $M(\A)$ and $M(\B)$.    Then for $\A, \B \in \Re^{m \times n \times m \times n}$, we have
\begin{equation} \label{e6.15}
\langle \A, \B \rangle = \langle M(\A), M(\B) \rangle,
\end{equation}
and
\begin{equation} \label{e6.16}
\|\A\|_2 = \|M(\A)\|_2.
\end{equation}

We first prove a proposition.

\begin{proposition} \label{p5.1}
Let $\A = (a_{i_1j_1i_2j_2}) \in \Re^{m \times n \times m \times n}$. Then
$$\| M(\A)\|_S \ge \| \A \|_S.$$
\end{proposition}
{\bf Proof}
Suppose that $\x \in \Re^m, \y \in \Re^n$. Let $\x \otimes \y$ be the Kronecker product of $\x$ and $\y$.  Then $\x \otimes \y \in \Re^{mn}$.  If $\|\x \|_2 = \|\y\|_2 =1$, then
$\|\x \otimes \y\|_2 = 1$.

By Theorem \ref{t01},  we have
$$\begin{aligned}
\|\A \|_S & = \max \left\{ \left| \langle \A, \x \circ \y \circ \uu \circ \vv \rangle \right| :  \|\x\|_2 = \|\y\|_2 = \|\uu\|_2 = \|\vv\|_2 = 1, \x, \uu \in \Re^m, \y, \vv \in \Re^n \right\}\\
& = \max \left\{ \left| \langle M(\A), (\x \otimes \y) \circ (\uu \otimes \vv) \rangle \right| :  \|\x\|_2 = \|\y\|_2 = \|\uu\|_2 = \|\vv\|_2 = 1, \x, \uu \in \Re^m, \y, \vv \in \Re^n \right\}\\
& \le  \max \left\{ \left| \langle M(\A), \z \circ  \w \rangle \right| :  \|\z\|_2 = \|\x\|_2  = 1, \z, \w \in \Re^{mn} \right\}\\
& = \|M(\A)\|_S.
\end{aligned}$$
This proves the proposition.
\qed

If $\A = (a_{i_1j_1i_2j_2}) \in BQ(m, n)$, then its matrix flattening $M(\A)$ is symmetric.   We now have the following theorem.

\begin{theorem}  \label{t6.2}
Suppose that $\A = (a_{i_1j_1i_2j_2}) \in BQ(m, n)$ and $M = M(\A)$ is its symmetric matrix flattening.  Then
\begin{equation} \label{e6.17}
\| M \|_* \le \|\A \|_* \le \min \{m, n \}\| M \|_*.
\end{equation}
In particular, if $\A$ is diagonal, we have
\begin{equation} \label{e6.18}
\|\A \|_* = \| M \|_* = \sum_{i=1}^m \sum_{j=1}^n |a_{ijij}|.
\end{equation}
\end{theorem}
{\bf Proof}  We first prove the first inequality of (\ref{e6.17}).   By (\ref{n2}), we have
$$\begin{aligned} \| \A \|_* & = \max \left\{ \left| \langle \A, \B \rangle \right| :  \|\B \|_S = 1, \B \in \Re^{m \times n \times m \times n} \right\}\\
& = \max \left\{ \left| \langle \A, \B \rangle \right| :  \|\B \|_S \le 1, \B \in \Re^{m \times n \times m \times n} \right\}\\
& = \max \left\{ \left| \langle M(\A), M(\B) \rangle \right| :  \|\B \|_S \le 1, \B \in \Re^{m \times n \times m \times n} \right\}\\
& \ge \max \left\{ \left| \langle M(\A), M(\B) \rangle \right| :  \|M(\B) \|_S \le 1, \B \in \Re^{m \times n \times m \times n} \right\}\\
& \ge \max \left\{ \left| \langle M(\A), B \rangle \right| :  \|B \|_S \le 1, B \in \Re^{mn \times mn} \right\}\\
& = \|M(\A)\|_*,
\end{aligned}$$
where the third equality is due to (\ref{e6.15}), the first inequality is by Proposition \ref{p5.1}, and the last equality is by the definition of the nuclear norm of a matrix.

We now prove the second inequality of (\ref{e6.17}).   Since $M(\A) \in \Re^{mn \times mn}$ is symmetric, we may assume that $M(\A)$ has an eigenvalue decomposition
$$M(\A) = \sum_{k=1}^{mn} \lambda_k \z^{(k)} \left(\z^{(k)}\right)^\top,$$
where $\z^{(k)} \in \Re^{mn}$ and $\| \z^{(k)} \|_2 = 1$, for $k = 1, \cdots, mn$.    For each $k$, $\z^{(k)}$ corresponds to an $m \times n$ matrix $\bar M(\z^{(k)})$.  Since
$\| \z^{(k)} \|_2 = 1$, we have  $\| \bar M(\z^{(k)}) \|_F = 1$, where $\| \cdot \|_F$ is the Frobenius norm.  Then \cite{Hu15}, we have
$$\|  \bar M(\z^{(k)}) \|_* \le \sqrt{\min\{ m, n\}}.$$
On the other hand,
$$\|  \bar M(\z^{(k)}) \|_* = \sum_{l=1}^{\min\{ m, n\}} |\sigma_{k, l}|,$$
where $\sigma_{k, l}$ for $l = 1, \cdots, \min \{ m, n \}$, are singular values of $\bar M(\z^{(k)})$ for $k = 1, \cdots, mn$.
Then $\bar M(\z^{(k)})$ has a singular value decomposition
$$\bar M(\z^{(k)}) = \sum_{l=1}^{\min\{ m, n\}} \sigma_{k, l} \x^{(k, l)}\left(\y^{(k, l)}\right)^\top,$$
where $\x^{(k, l)} \in \Re^m$, $\|\x^{(k, l)} \|_2 = 1$  and $\y^{(k, l)} \in \Re^n$, $\|\y^{(k, l)} \|_2 = 1$, for $k = 1, \cdots, mn$ and $l = 1, \cdots, \min \{ m, n \}$.
This implies
$$\sum_{l=1}^{\min\{ m, n\}} |\sigma_{k, l}| \le  \sqrt{\min\{ m, n\}},$$
for $k = 1, \cdots, mn$.   Then we have
$$\begin{aligned} \A & = \sum_{k=1}^{mn} \lambda_k \bar M(\z^{(k)}) \circ \bar M(\z^{(k)})\\
& = \sum_{k=1}^{mn} \lambda_k \left(\sum_{l=1}^{\min\{ m, n\}} \sigma_{k, l} \x^{(k, l)} \circ \y^{(k, l)}\right) \circ \left(\sum_{l=1}^{\min\{ m, n\}} \sigma_{k, l} \x^{(k, l)} \circ \y^{(k, l)}\right)\\
& = \sum_{k=1}^{mn} \lambda_k \sum_{l, s=1}^{\min\{ m, n\}} \sigma_{k, l}\sigma_{k, s} \x^{(k, l)} \circ \y^{(k, l)} \circ \x^{(k, s)} \circ \y^{(k, s)}.
\end{aligned}$$
Thus,
$$\|\A\|_* \le \sum_{k=1}^{mn} \sum_{l, s=1}^{\min\{ m, n\}} |\lambda_k \sigma_{k, l}\sigma_{k, s}| \le \sum_{k=1}^{mn} |\lambda_k|\min\{ m, n\} = \min\{ m, n\}\| M(\A)\|_*. $$

Finally, assume that $\A$ is diagonal.  Then
$$\A = \sum_{i=1}^m \sum_{j=1}^n a_{ijij} \ee^{(i)} \circ \bar \ee^{(j)} \circ \ee^{(i)} \circ \bar \ee^{(j)},$$
where $\ee^{(i)}$ for $i = 1, \cdots m$ are the unit vectors in $\Re^m$, while $\bar \ee^{(j)}$ for $j = 1, \cdots n$ are the unit vectors in $\Re^n$.
This implies that
$$\|\A\|_* \le \sum_{i=1}^m \sum_{j=1}^n |a_{ijij}| = \|M(\A)\|_*.$$
Then, by (\ref{e6.17}), we have (\ref{e6.18}).

This proves the theorem.
\qed

Comparing this theorem with Theorem 5.2 of \cite{Hu15}, our theorem is somewhat stronger.

This theorem says that the equality in the first inequality of (\ref{e6.17}) may hold.  How about the second inequality of (\ref{e6.17})?

\begin{corollary} \label{c5.3}
$$\|\II\|_* = mn$$.
\end{corollary}

\section{Norms of Tensor Products of Biquadratic Tensors}

We may define products of two biquadratic tensors.  Let $\A = (a_{i_1j_1i_2j_2}), B = (b_{i_1j_1i_2j_2}) \in BQ(m, n)$, then we have
$\C = (c_{i_1j_1i_2j_2}) : = \A \B \in BQ(m, n)$, defined by
$$c_{i_1j_1i_2j_2} = \sum_{i_3=1}^m \sum_{j_3=1}^n a_{i_1j_1i_3j_3}b_{i_3j_3i_2j_2},$$
for $i_1, i_2 = 1, \cdots, m$ and $j_1, j_2 = 1, \cdots, n$.

   Then, for any $\A \in BQ(m, n)$,
$$\A\II = \II\A = \A.$$
If $\A, \B \in BQ(m, n)$ and $\A \B = \II$, then we also have $\B \A = \II$ and we denote $\A^{-1} = \B$.

We have the following proposition.

\begin{proposition}   \label{p6.1}
For any $\A, \B \in BQ(m, n)$, we have
$$\| \A \B \|_* \le \|\A\|_*\|\B\|_*.$$
\end{proposition}

This proposition may be proved directly.  It may also be regarded a special case of Theorem 2.1 of \cite{QHZ19}.  Hence, we do not prove it here.

By Proposition \ref{p6.1} and Corollary \ref{c5.3}, we have the following proposition.

\begin{proposition}
Suppose that $\A \in BQ(m, n)$ is invertible.  Then we have
$$\|\A\|_* \|\A^{-1}\|_* \ge mn.$$
\end{proposition}

In general, for $\A, \B \in BQ(m, n)$, we may not have
$$\|\A \B\|_S \le \|\A\|_S \|\B\|_S.$$
See Example 4.1 of \cite{QHZ19}.  On the other hand, for any $\A, \B \in BQ(m, n)$, by Theorem 4.2 of \cite{QHZ19}, we have
$$\|\A \B\|_S \le \|\A\|_* \|\B\|_S.$$
By definition, it is easy to see that
$$\|\II\|_S = 1.$$
From these, we have the following proposition.

\begin{proposition}
Suppose that $\A \in BQ(m, n)$ is invertible.  Then we have
$$\|\A\|_* \|\A^{-1}\|_S \ge 1.$$
\end{proposition}

\section{Final Remarks}

Viewing the importance and the special structure properties of biquadratic tensors, we hope that we may explore more at this direction.

Our study can be extended to the field of complex numbers without difficulties.

Our study can also be extended to bisymmetric tensors.   A former definition for bisymmetric tensors are as follows.

\begin{definition}
Let $p$ be a positive integer.  Let $\Re^{n_1 \cdots \times n_p \times n_1 \cdots \times n_p}$ be the space of $2p$th order tensors of dimension $n_1 \cdots \times n_p \times n_1 \cdots \times n_p$.   Let $\A = (a_{i_1\cdots i_pj_1\cdots j_p}) \in \Re^{n_1 \cdots \times n_p \times n_1 \cdots \times n_p}$.  The tensor $\A$ is called bisymmetric if for all $i_k, j_k = 1, \cdots, n_k$, $k = 1, \cdots, p$, we have
$$a_{i_1\cdots i_pj_1\cdots j_p} = a_{i_1\cdots i_{k-1}j_ki_{k+1}\cdots i_pj_1\cdots j_{k-1}i_kj_{k+1}\cdots j_p}.$$
\end{definition}
Then for $p=1$, we have symmetric matrices, and for $p=2$, we have biquadratic tensors.

\bigskip

{\bf Acknowledgment}   The authors are thankful to
Bo Jiang for the discussion on Theorem 4.4, and to Chen Ling for his comments.

\end{document}